\documentclass{qam-l}

\usepackage{amssymb,amsmath} 

\allowdisplaybreaks[1]

\renewcommand{\appendix}{
  \setcounter{section}{0}\renewcommand{\thesection}{\Alph{section}}
  \section*{Appendix} 
}

\def\Appendix#1{
  \setcounter{equation}{0}
  \renewcommand{\theequation}{\thesection.\arabic{equation}}
  \section{#1}
}

\newtheorem{theorem}{Theorem}[section]
\newtheorem{lemma}[theorem]{Lemma}
\theoremstyle{definition}

\theoremstyle{remark}

\numberwithin{equation}{section}

\newcommand{\rf}[1]{(\ref{#1})}

\def\bxt{\! \boxtimes \! }

\def\bd#1{\mbox{\boldmath$\displaystyle\mathbf{#1}$} }
\def\tens#1{\mathbb{\,#1}}	   
\def\tr{\operatorname{tr}} 
\def\dd{\operatorname{d}}

\usepackage[pdftitle={tensors, gradients,continuum mechanics},
              pdfauthor={Andrew N. Norris},
              pdfkeywords={tensors, gradients,continuum mechanics}, 
              letterpaper=true,
              pdfpagetransition=Dissolve]{hyperref}
\begin{document} 

\pagestyle{myheadings}\markright{{\sc  Norris }  ~~~~~~\today}

\title{Higher  derivatives and the inverse derivative of a tensor-valued function of a tensor}

\author{Andrew N. Norris}
\address{Mechanical and Aerospace Engineering, 	Rutgers University, Piscataway NJ 08854-8058, USA}
\email{norris@rutgers.edu}
	\date{}

\begin{abstract}

The  $n^{th}$ derivative of a tensor valued function of a tensor is defined by a finite number of coefficients each with closed form expression. 


\end{abstract}
\maketitle

\section{Introduction}

We consider tensor functions on symmetric second order tensors, $S$ym$\rightarrow S$ym, defined by a  scalar function $f(x)$ of a single variable according to 
\begin{equation} \label{-0}
f({\bd A}) = \sum\limits_{i=1}^d f(\alpha_i) {\bd A}_i,
\end{equation}
where $d$ is the eigen-index of ${\bd A}$. 
The tensor  ${\bd A}\in S$ym is   arbitrary  with spectral decomposition
\begin{equation} \label{-71}
{\bd A}  = \sum\limits_{i=1}^d \alpha_i {\bd A}_i ,  
\qquad  
{\bd A}_i{\bd A}_j = \begin{cases}
{\bd A}_i & i=j,\\
0 , & i\ne j,  
\end{cases}
\qquad 
\sum\limits_{i=1}^d  {\bd A}_i = {\bd I} .  
\end{equation}
The specific case of tensors acting on 3-dimensional vectors, $d\le 3$, is discussed in this paper, although the results can be readily generalized. 

Derivatives of $f( {\bd A}) $ are defined by the expansion
\begin{align} \label{-81}
 f( {\bd A} + {\bd X}) =& f({\bd A})  + \nabla f({\bd A}) {\bd X}
  + \frac12 \nabla^{(2)} f({\bd A}) : {\bd X}{\bd X}
   + \frac1{3!} \nabla^{(3)} f({\bd A}) : {\bd X}{\bd X}{\bd X}
   \nonumber \\ & \quad 
+ \ldots + 
 \frac1{n!} \nabla^{(n)} f({\bd A}) : 
  \underbrace{ {\bd X}{\bd X}\ldots {\bd X}}_{n \text{ terms}}
 +\ldots
 \end{align}
 The $n^{th}$ derivative $\nabla^{(n)} f({\bd A})$ is a tensor of order $2(n+1)$ which contracts $n$-times  with the  second order tensor ${\bd X}$ to produce a second order tensor. 
The first derivative, or gradient, is  
\begin{equation} \label{45}
 \nabla f ({\bd A})   = 
 \sum\limits_{i,j=1}^d \, f^{(1)}_{ij} 
 {\bd A}_i \bxt {\bd A}_j, 
 \qquad
 f^{(1)}_{ij}    = \begin{cases} 
 f'(\alpha_i),  &  i=j ,  \\ 
  & \\
 \frac{ f(\alpha_i) - f(\alpha_j)}{ \alpha_i - \alpha_j}, & i\ne j\text ,
 \end{cases} 
 \end{equation}
 where $\bxt$ denotes the outer tensor product, defined  in Section \ref{sec2}. 
 The identity \rf{45} is well known and has appeared in various formats.  The first explicit presentation I am aware of is due to Ogden \cite{Ogden84} who defines a   fourth order tensor 
 $\bd{\mathcal L}^1 = \partial f ( {\bd A}) /  \partial  {\bd A} $ (with slight change in notation). 
Ogden  gives the components ${\mathcal L}^1_{ijkl}$ in terms of the eigenvectors  ${\bd a}_i$ of ${\bd A}$, i.e. 
 $\bd{\mathcal L}^1 = {\mathcal L}^1_{ijkl} {\bd a}_i\otimes {\bd a}_j\otimes {\bd a}_k\otimes {\bd a}_l$. 
 These coefficients are related to those of \rf{45} by the 
fact that  ${\bd A}_i = {\bd a}_i\otimes{\bd a}_i$ (no sum) when $d=3$ and hence 
${\mathcal L}^1_{iiii} = f^{(1)}_{ii} $, ${\mathcal L}^1_{ijij} = f^{(1)}_{ij} $.  
The fundamental result \rf{45} was also derived by  Carlson and Hoger \cite{Carlson86}, although the present notation  is based on  \cite{Xiao98}. 

Ogden \cite{Ogden84} (Section 3.4) also presented the second derivative.  In the present notation it is
\begin{align} \label{47}
 \frac12 \nabla^{(2)} f ({\bd A})   &= 
 \sum\limits_{i,j,k=1}^d \, f^{(2)}_{ijk} 
 {\bd A}_i \bxt {\bd A}_j\bxt {\bd A}_k, 
 \\
 f^{(2)}_{iii}    &=  \frac12 f''(\alpha_i),  
  \nonumber \\
  f^{(2)}_{iij}    &= 
  \frac{ f(\alpha_j) - f(\alpha_i) - (\alpha_j - \alpha_i) f'(\alpha_i)}
{(\alpha_j - \alpha_i )^2}, 
 \nonumber \\
  f^{(2)}_{ijk}    &= 
  \frac{ [ f(\alpha_j) - f(\alpha_i) ]( \alpha_i+\alpha_j-2 \alpha_k)
  - [ f(\alpha_i) + f(\alpha_j) -2 f(\alpha_k) ] ( \alpha_j- \alpha_i)}
  {2( \alpha_i- \alpha_j)( \alpha_j- \alpha_k)( \alpha_k- \alpha_i)}.
  \nonumber 
 \end{align}
 Ogden's  expressions \rf{45} and \rf{47} are special cases of the more general formulae for 
derivatives of isotropic tensor functions derived by Chadwick and Ogden \cite{Chadwick71}. 

 It is interesting to note that the tensor gradient function (first derivative) involves finite differences of the function $f$ as well as its derivative evaluated at the eigenvalues of $\bd A$. 
Similarly, the second derivative contains second order finite differences.  The general result derived here shows that the coefficients for the $n^{th}$ derivative are related to  an interpolating polynomial.  
 
Our main result is the following: 
\begin{theorem}\label{thm1}
The $n^{th}$  derivative of the tensor function $f({\bd A})$ is given by 
\begin{equation} \label{54}
\frac1{n!} \nabla^{(n)} f({\bd A}) =  
\sum\limits_{  i_1,i_2,\ldots ,  i_{n+1}=1  }^d  \,
 f^{(n)}_{i_1,i_2,\ldots ,  i_{n+1}} \, 
 {\bd A}_{i_1} \bxt {\bd A}_{i_2} \bxt \ldots \bxt {\bd A}_{i_{n+1}} .
   \end{equation}
 The coefficients $f^{(n)}_{i_1,i_2,\ldots ,  i_{n+1}}$ are unaltered under all permutations of the $n+1$ indices.  The $(n+2)(n+3)/2$  distinct coefficients can be classified into\footnote{ $\lfloor x \rfloor $ is the floor function.} 
 $\lfloor \frac{ (n+4)^2 + 4}{12} \rfloor $ 
  expressions $f_{ijk }^{\nu_i ,\nu_j ,\nu_k }$ where 
 $\nu_i, \nu_j, \nu_k$, are the number of   occurrences of distinct indices $i,j,k$ respectively, with $i\ne j\ne k\ne i$ and  $\nu_i+ \nu_j+ \nu_k=n+1$.  The coefficient is
 \begin{equation} \label{305}
f_{ijk }^{\nu_i ,\nu_j ,\nu_k } = \sum\limits_{\stackrel{l=i,j,k}{\nu_l > 0}} \frac{1}{(\nu_l-1)!}
 \frac{\dd^{\nu_l-1} ~}{\dd x^{\nu_l-1}} 
  \frac{f(x)}{ \prod\limits_{\stackrel{m=i,j,k}{m\ne l}} (x-\alpha_m)^{\nu_m} } 
\bigg|_{x=\alpha_l}  .
 \end{equation}
  Alternatively, the coefficient can be found from  the unique interpolating polynomial $P(x)$ of degree $n$ that fits the data  at the three points  $x=  \{\alpha_i, \alpha_j, \alpha_k\}$ defined by  the $n+1$ values 
 $f^{(I)}(\alpha_i)$,  $f^{(J)}(\alpha_j)$, 
 and  $f^{(K)}(\alpha_k)$ for  $ 0\le I \le \nu_i-1 , 0\le J \le \nu_j-1, 0\le K \le \nu_k-1 $
 where
 $f^{{(l)}}(x)$ is the $l^{th}$ derivative. 
Let $P(x) = p_n x^n+p_{n-1}x^{n-1} + \ldots + p_0$, then 
\begin{equation} \label{090}
f_{ijk }^{\nu_i ,\nu_j ,\nu_k } = p_n.
 \end{equation}
 \end{theorem}

The first few expressions for the coefficients are 
\begin{subequations}\label{94}
\begin{align} \label{1a}
f_{ijk}^{0,0,n+1} &= \frac{1}{n!} f^{(n)}(\alpha_k), 
\\
f_{ijk}^{0,1,n}~~~&= \frac{1}{(\alpha_j - \alpha_k)^{n-1}}
\bigg[
f(\alpha_j) -  \sum\limits_{l=0}^{n-1}
\frac{1 }{l!}(\alpha_j - \alpha_k)^l f^{(l)}(\alpha_k)
\bigg],
\label{1b}
\\
f_{ijk}^{0,2,n-1}&= \frac{1}{(\alpha_j - \alpha_k)^{n-1}}
\bigg[ \sum\limits_{l=0}^{n-1}\frac{(n-l)} {l!}  (\alpha_j - \alpha_k)^l f^{(l)} (\alpha_k)
+ (\alpha_j - \alpha_k) f'(\alpha_j)- n f(\alpha_j)
\bigg],
\label{1c}
\\
f_{ijk}^{1,1,n-1} &= \frac{1}{(\alpha_i - \alpha_j) }
\bigg\{
\frac{ f(\alpha_i)}{ (\alpha_i - \alpha_k)^{n-1} } -
\frac{ f(\alpha_j)}{ (\alpha_j - \alpha_k)^{n-1} }
\nonumber \\
& \qquad\qquad\qquad\qquad 
- \sum\limits_{l=0}^{n-2}
\frac{1 }{l!}  f^{(l)}(\alpha_k) \bigg[
\frac{ 1}{ (\alpha_i - \alpha_k)^{n-1-l} }-
\frac{ 1}{ (\alpha_j - \alpha_k)^{n-1-l} }
\bigg]
\bigg\},
\label{1d}
\\
f_{ijk}^{1,2,n-2} &= \frac{1}{(\alpha_i - \alpha_j)^2 }
\bigg\{
\frac{ f(\alpha_i)}{ (\alpha_i - \alpha_k)^{n-2} } -
\frac{ f(\alpha_j)}{ (\alpha_j - \alpha_k)^{n-2} }
\big[ 1+ (n-2)\frac{(\alpha_j - \alpha_i)}{(\alpha_j - \alpha_k)}\big]
\nonumber \\
& \qquad\qquad\qquad + \frac{ (\alpha_j - \alpha_i)  }{ (\alpha_j - \alpha_k)^{n-2} }f'(\alpha_j)
- \sum\limits_{l=0}^{n-3}\frac{1 }{l!}  f^{(l)}(\alpha_k) \bigg[
\frac{ 1}{ (\alpha_i - \alpha_k)^{n-2-l} } 
\nonumber \\
& \qquad\qquad\qquad -
\frac{ 1}{ (\alpha_j - \alpha_k)^{n-2-l} }
\big[ 1+ (n-l-2)\frac{(\alpha_j - \alpha_i)}{(\alpha_j - \alpha_k)}\big]
\bigg]
\bigg\}.
\label{1e}
\end{align}
\end{subequations}

These  are sufficient to determine  all derivatives up to and including the fourth.  Thus, 
the gradient tensor $(n=1)$ requires the two expressions $f_{ijk}^{0,0,2}$ and $f_{ijk}^{0,1,1}$ evident from \rf{45}; the second derivative $(n=2)$ requires $f_{ijk}^{0,0,3}$, 
$f_{ijk}^{0,1,2}$ and $f_{ijk}^{1,1,1}$, which may be read off from \rf{47}; the third  derivative $(n=3)$ involves four distinct formulas for $f_{ijk}^{0,0,4}$, 
$f_{ijk}^{0,1,3}$, $f_{ijk}^{0,2,2}$ and $f_{ijk}^{1,1,2}$; and the   fourth  derivative $(n=4)$ requires the five expressions  $f_{ijk}^{0,0,5}$, $f_{ijk}^{0,1,4}$, $f_{ijk}^{0,2,3}$, 
$f_{ijk}^{1,1,3}$ and $f_{ijk}^{1,2,2}$.  
Note that \rf{1b} and \rf{1c} reduce to \rf{1a} in the limit as $\alpha_i\rightarrow \alpha_k$.  Similarly, \rf{1d} and \rf{1e} reduce to \rf{1c}  in the limit.

The main objective of this paper is to prove Theorem \ref{thm1}.  Some new results concerning the properties of the fourth order gradient tensor $\nabla f$ are  also presented and the inverse fourth order tensor $\nabla^{-1} f$ is introduced.  Both the gradient and its inverse are discussed with application to strain measure functions \cite{Hill78,Scheidler91}. 
The proof of  Theorem \ref{thm1} begins with   a new derivation of the well known expression for the gradient $\nabla f$. 
The essential structure of the second and higher order derivatives is shown to depend  
on a general algebraic identity.  This  identity also reveals the appearance of the characteristic finite difference terms.  The proof is completed by making  connections with contour integrals and with interpolation polynomials.  The results are presented in terms of Kronecker products of tensors which makes the expressions more transparent.  

The paper is laid out as follows. Notation is introduced in Section \ref{sec2} followed by 
the derivation of Theorem \ref{thm1}  in Section \ref{sec3}.    The inverse gradient tensor is introduced and its properties discussed in Section \ref{sec4}.  

\section{Notation and preliminaries }\label{sec2}

We consider second order tensors acting on vectors in 
a three dimensional inner product space,  
 ${\bd x} \rightarrow {\bd A}{\bd x}$ with  transpose ${\bd A}^t$ such that 
 ${\bd y}\cdot{\bd A}{\bd x} = {\bd x}\cdot{\bd A}^t{\bd y}$.  Spaces of symmetric and skew-symmetric tensors are distinguished, Lin = Sym $\oplus$ Skw 
where   ${\bd A}\in $ Sym (Skw) iff 
${\bd A}^t={\bd A}$ (${\bd A}^t=-{\bd A}$).  
Products ${\bd A} {\bd B}\in$Lin are defined by 
${\bd y}\cdot{\bd A} {\bd B}{\bd x} =   ({\bd A}^t{\bd y}) \cdot {\bd B}{\bd x}$. 

Psym  is the space of positive definite second order tensors.
Functions of a symmetric tensor 
can be phrased in terms of its spectral decomposition \rf{-0} where ${\bd A}_i\in$Psym and the distinct eigenvalues $\alpha_i$, $i=1\ldots , d\le 3$ are real numbers. 
The single function $f({\bd A})$  is a special case of isotropic tensor functions of the form 
$
{\bd T}({\bd A}) = \sum\limits_{i=1}^3 f_i(\alpha_1,\alpha_2,\alpha_3 ) {\bd A}_i 
$,
involving three functions of three variables.  Chadwick and Ogden \cite{Chadwick71} derived first and second derivatives for   more general situation (see \cite{Chen04} for a recent overview).   Our interest here is strictly limited to one function $f$ as in \rf{-0}. 

${\tens L}$in is the space of fourth order tensors acting on Lin. 
The square tensor product or Kronecker product ${\bd A}\bxt {\bd B}$, Lin$\times$Lin$\rightarrow {\tens L}$in, is defined in the usual manner as  \cite{Rosati00}
\begin{equation} \label{51}
( {\bd A}\bxt {\bd B}) {\bd X} = {\bd A}{\bd X}{\bd B}^t, 
\qquad \forall \,  {\bd X}\in \text{ Lin}. 
\end{equation}
The following property of $\bxt$ will be used extensively, 
\begin{equation} \label{53}
( {\bd A}\bxt {\bd B}) ( {\bd X}\bxt {\bd Y})    = 
( {\bd A}{\bd X}) \bxt ({\bd B}{\bd Y}) .
\end{equation}
The generalization of \rf{51} is 
\begin{equation} \label{04}
\big(
\underbrace{{\bd A}\bxt{\bd B}\ldots  \bxt\ldots  {\bd C}}_{n+1}\big)   : 
\underbrace{ {\bd X}{\bd Y}\ldots {\bd Z} }_{n}
= {\bd A}{\bd X}{\bd B}^t {\bd Y}\ldots  {\bd Z} {\bd C}^t. 
\end{equation}
Basic properties of the Kronecker product can be found in e.g. \cite{Kintzel06}.

The gradient of a tensor function $f({\bd A})$  is a fourth order tensor $\nabla f \in \tens{L}$in defined by 
\begin{equation} \label{44}
 \nabla f ({\bd A}) \, {\bd X} = 
 \lim\limits_{\epsilon \rightarrow 0} \frac1{\epsilon }\big[ 
 f ({\bd A} +\epsilon {\bd X}) - f ({\bd A})
 \big]. 
 \end{equation}
Higher order derivatives are defined recursively in accordance with \rf{-81}, 
\begin{equation} \label{144}
 \nabla^{(n)} f ({\bd A}) \, {\bd X} = 
 \lim\limits_{\epsilon \rightarrow 0} \frac1{\epsilon }\big[ 
 \nabla^{(n-1)}  f ({\bd A} +\epsilon {\bd X}) - \nabla^{(n-1)} f ({\bd A})
 \big], \qquad n\ge 2. 
 \end{equation}

We note some properties of the tensor gradient. First,  the tensor gradient of the product of two functions is   
\begin{equation} \label{558}
\nabla [ f({\bd A}) g({\bd A}) ] = \big({\bd I} \bxt g({\bd A})\big) \nabla  f({\bd A})
+ \big({f({\bd A})\bxt \bd I}  \big) \nabla  g({\bd A}). 
\end{equation}
For example, let $g= 1/f \equiv f^{-1}$, then using   $\nabla [ f({\bd A}) g({\bd A}) ] = 0$
yields
\begin{equation} \label{559}
\nabla   f^{-1}({\bd A}) = - \big[f^{-1}({\bd A}) \bxt f^{-1}({\bd A})\big] \nabla  f({\bd A}). 
\end{equation}
The composition of two functions is $f\circ g({\bd A}) = f\big( g({\bd A})\big) $ has gradient
\begin{equation} \label{557}
\nabla  f\circ g({\bd A}) = 
\nabla  f\big( g({\bd A})\big)  \, \nabla g({\bd A}). 
\end{equation}
Only symmetric tensor arguments are considered in this paper.  
The extension of the gradient tensor to non-symmetric tensor arguments is discussed by Itskov  \cite{Itskov02b}. 

\section{Proof of Theorem \ref{thm1} }\label{sec3}

The proof has several  stages, starting with consideration of $f$ in the form of a power series.  This allows us to identify certain properties of the  derivatives that are subsequently  generalized  to arbitrary $f$.  Implicit in this approach is the idea  that a function of a single variable can be approximated to any desired level of precision by a polynomial.   The use of polynomials allows us to see the structure of the higher order derivatives, not only the first or second.   

\subsection{A power series}
For the moment assume 
\begin{equation} \label{811}
f(x) =   x^m ,  
\end{equation}
then $f( {\bd A} + {\bd X})$ is 
\begin{align} \label{01}
( {\bd A} + {\bd X})^m =& ( {\bd A} + {\bd X})( {\bd A} + {\bd X})(\ldots ) ( {\bd A} + {\bd X})
\nonumber \\
=&{\bd A}^m + \big( {\bd A}^{m-1} {\bd X} + {\bd A}^{m-2} {\bd X}{\bd A}
+ \ldots  + {\bd X}{\bd A}^{m-1}\big) 
\nonumber \\
 & + \big( {\bd A}^{m-2} {\bd X}^2 + {\bd A}^{m-3} {\bd X}{\bd A}{\bd X}
+ \ldots  + {\bd X}^2{\bd A}^{m-2}\big) + \ldots . 
\end{align}
The O$({\bd X})$ and O$({\bd X}^2)$ contributions  contain $m$ and $m(m-1)/2$ separate terms, respectively. The term of O$({\bd X}^k)$ has $\binom{m}{k}$ elements. 
A more suggestive and useful form of  \rf{01} is apparent by rewriting it as
\begin{align} \label{03}
( {\bd A} + {\bd X})^m - {\bd A}^m =& \sum\limits_{k=1}^{m}
 {\bd A}^{m-k} {\bd X} {\bd A}^{k-1}  
  + \sum\limits_{k=1}^{m-1}
 {\bd A}^{m-1-k} {\bd X} \sum\limits_{l=1}^{k} {\bd A}^{l-1} {\bd X} {\bd A}^{k-l}  
 +\ldots 
 \nonumber \\
 =& \sum\limits_{k=1}^{m} \big( {\bd A}^{m-k} \bxt {\bd A}^{k-1} \big) {\bd X}
 + \sum\limits_{k=1}^{m-1} \sum\limits_{l=1}^{k}  
 {\bd A}^{m-1-k} \bxt  {\bd A}^{l-1} \bxt {\bd A}^{k-l}   : {\bd X}{\bd X}
  \nonumber \\
  & + \sum\limits_{ \stackrel{i_i,i_2, \ldots , i_{n+1}=0} {i_i+  \ldots + i_{n+1}=m-n}}^{m-n} 
x_1^{i_i}\, x_2^{i_2}  \, \ldots \, x_{n}^{i_{n}} \, x_{n+1}^{i_{n+1}}
 {\bd A}^{i_i}\bxt {\bd A}^{i_2} \bxt \ldots \bxt {\bd A}^{i_{n}} \bxt {\bd A}^{i_{n+1}}
 \nonumber \\
  & \qquad \qquad \qquad \qquad : \underbrace{{\bd X}{\bd X} \ldots {\bd X} }_{n} 
 +\ldots .
\end{align}
This permits us to consider the derivatives in sequence.  Our objective is to understand the 
term of O$({\bd X}^n)$, but it is easier and more instructive to begin with the lowest, $n=1$.

\subsection{The first derivative} 
Let us first focus on the  O$({\bd X})$ term.   Consider the identity 
\begin{align} \label{05}
({\bd A}   \bxt {\bd I} - {\bd I} \bxt {\bd A}) \, 
\sum\limits_{k=1}^{m}   {\bd A}^{m-k} \bxt {\bd A}^{k-1} 
  &= \sum\limits_{k=1}^{m} ( {\bd A}^{m+1-k} \bxt {\bd A}^{k-1} - 
  {\bd A}^{m-k} \bxt {\bd A}^{k} )
  \nonumber \\
 &=  {\bd A}^m \bxt {\bd I} - {\bd I} \bxt {\bd A}^m.
\end{align}
Then
\begin{equation} \label{812}
({\bd A}   \bxt {\bd I} - {\bd I} \bxt {\bd A}) \, 
\big[
( {\bd A} + {\bd X})^m - {\bd A}^m \big] = ({\bd A}^m \bxt {\bd I} - {\bd I} \bxt {\bd A}^m) \, {\bd X} 
 + \text{O}({\bd X}^2) , 
\end{equation}
which may be generalized to power series functions 
\begin{equation} \label{810}
f(x) = \sum_{m=0}^M c_m x^m, 
\end{equation}
as 
\begin{equation} \label{813}
({\bd A} \bxt {\bd I} - {\bd I} \bxt {\bd A}) \, 
\big[
f( {\bd A} + {\bd X}) - f({\bd A}) \big] 
= \big[ f( {\bd A})\bxt {\bd I} - {\bd I} \bxt f( {\bd A}) \big] {\bd X} 
 + \text{O}({\bd X}^2).
\end{equation}
Equation \rf{813} implies 
\begin{lemma}\label{lem1}
The tensor gradient of  any differentiable function satisfies
 \begin{equation} \label{8131}
({\bd A} \bxt {\bd I} - {\bd I} \bxt {\bd A}) \, 
\nabla f({\bd A}) 
=   f( {\bd A})\bxt {\bd I} - {\bd I} \bxt f( {\bd A}) .
\end{equation}
\end{lemma}

We will return to the form of this equation later when  the inverse gradient tensor is discussed. At this stage we note  an isomorphism between the  tensorial identity \rf{05} and an algebraic equation in two variables. 
Let $x$ and $y$ represent the first and second occurrence of $\bd A$ in 
${\bd A}^{m-k} \, \bxt {\bd A}^{k-1} \rightarrow x^{m-k}y ^{k-1}$.  Then eq. \rf{05}  is equivalent to  
\begin{align} \label{06}
(x-y) \sum\limits_{k=1}^{m}
 x^{m-k} y^{k-1}
&= 
(x-y) (x^{m-1} + x^{m-2}y+ x^{m-3}y^2 + \ldots  +  y^{m-1})
\nonumber \\
&=  x^m-y^m. 
\end{align}
This  is the well known factorization of $x^m-y^m$.    The 
analogous relation for the second derivative of $f$ is  derived next. 

\subsection{The second derivative} 

The algebraic  equality analogous to \rf{06} for  
the O$({\bd X}^2)$ terms involves three variables, say $x$, $y$ and $z$, and is 
\begin{equation} \label{07}
(x-y)(y-z)(z-x) 
\sum\limits_{k=1}^{m-1}
  \sum\limits_{l=1}^{k} x^{m-1-k} y^{l-1} z^{k-l}
 = (z-y)x^m +(x-z)y^m +(y-x)z^m. 
 \end{equation} 
As proof of this statement note that the tensorial equivalent of \rf{07} is  
 \begin{align} \label{8132}
&({\bd A} \bxt {\bd I}\bxt{\bd I} - {\bd I} \bxt {\bd A}\bxt{\bd I}) 
({\bd I}\bxt {\bd A} \bxt {\bd I} - {\bd I} \bxt{\bd I}\bxt {\bd A}) 
({\bd I}\bxt {\bd I} \bxt {\bd A} - {\bd A} \bxt{\bd I}\bxt {\bd I}) 
\frac12 \nabla^{(2)} f({\bd A}) 
\nonumber \\
&  = f({\bd A})\bxt ({\bd I} \bxt {\bd A} - {\bd A} \bxt{\bd I}) + 
\big( {\bd A} \bxt f({\bd A})\bxt{\bd I} - {\bd I} \bxt f({\bd A})\bxt{\bd A}
\big) +
({\bd I} \bxt {\bd A} - {\bd A} \bxt{\bd I})\bxt f({\bd A}).
\end{align}
It may be checked, using \rf{03}, that this indeed holds.   

The algebraic identities   for the $n=1$ and $n=2$ derivatives are, respectively, 
\begin{subequations}\label{0800}
 \begin{align} \label{08}
 \sum\limits_{k=1}^{m}
 x^{m-k} y^{k-1} &= \frac{x^m}{x-y} + \frac{y^m}{y-x},
 \\
\sum\limits_{k=1}^{m-1}
  \sum\limits_{l=1}^{k} x^{m-1-k} y^{l-1} z^{k-l}
 &= \frac{x^m}{(x-y)(x-z)} +\frac{y^m}{(y-x)(y-z)} +\frac{z^m}{(z-x)(z-y)} .  
 \end{align}
 \end{subequations}
This form of the identity will prove more useful in the general case of the $n^{th}$ derivative, considered next. 
 
 \subsection{Higher derivatives}
 
 Referring to eq. \rf{03} it is clear that the term of O$({\bd X}^n)$ involves  $\binom{m}{k}$ products  of the form 
  \begin{equation} \label{000}
 {\bd A}^{i_i}\, {\bd X}\, {\bd A}^{i_2}  \, {\bd X}\, \ldots \, {\bd X}\, {\bd A}^{i_{n}} \, {\bd X}\, {\bd A}^{i_{n+1}}, 
 \qquad \text{ with }i_i+i_2+ \ldots + i_{n+1}=m-n.
 \end{equation}
 The analogous algebraic quantity  has a closed form expression similar to that found for the first and second derivatives in \rf{0800}.  It is given by the following general  identity  among $n+1$ independent variables: 
 \begin{lemma}\label{lem2}
 \begin{equation} \label{092}
\sum\limits_{ \stackrel{i_i,i_2, \ldots , i_{n+1}=0} {i_i+i_2+ \ldots + i_{n+1}=m-n}}^{m-n} 
x_1^{i_i}\, x_2^{i_2}  \, \ldots \, x_{n}^{i_{n}} \, x_{n+1}^{i_{n+1}}
 = \sum\limits_{i=1}^{n+1} 
 \frac{x_i^m}{\prod\limits_{\stackrel{j=1}{j\ne i}}^{n+1}  (x_i-x_j)}   .  
 \end{equation} 
 \end{lemma}
  As proof\footnote{Thanks to Doron Zeilberger for the proof.}, consider  the identity
 \begin{equation} \label{3-}
 \frac{1}{ (1-x_1t)(1-x_2t)(\ldots)(1-x_{n+1}t)}
 = \frac{1}{t^n} \sum\limits_{i=1}^{n+1}   \frac{1}{ (1-x_i t) }
  \frac{1}{\prod\limits_{\stackrel{j=1}{j\ne i}}^{n+1}  (x_i-x_j)}   .  
 \end{equation}
 The right hand side is simply the  expansion of the left in partial fractions.  Equation \rf{092} 
  follows by comparing the coefficients of $t^{m-n}$ on either side of \rf{3-}.


Motivated by the expansion \rf{03}, consider the following ansatz for the derivatives,
\begin{subequations}
\begin{align} \label{4a}
\nabla f({\bd A}) &= \sum\limits_{ i,j=1}^d  \, f_{ij}^{(1)}\, 
 {\bd A}_i \bxt {\bd A}_j , 
\\
\frac12 \nabla^{(2)} f({\bd A}) &= \sum\limits_{ i,j,k=1}^d  \, f_{ijk}^{(2)}\, 
 {\bd A}_i \bxt {\bd A}_j\bxt {\bd A}_k , 
 \label{4b}
 \\  
\frac1{n!} \nabla^{(n)} f({\bd A}) &= \sum\limits_{ \stackrel{ i_1,i_2,\ldots ,}{\ldots , i_n ,   i_{n+1}=1} }^d  \,
 f_{i_1,i_2,\ldots ,  i_{n+1}}^{(n)} \, 
 {\bd A}_{i_1} \bxt {\bd A}_{i_2} \bxt \ldots \bxt {\bd A}_{i_{n+1}} .
  \label{4c}
\end{align}
\end{subequations}
The key feature is the use of   the spectral basis of $\bd A$.  

The tensorial equations \rf{8131} and \rf{8132}, or their algebraic versions \rf{06} and \rf{07},  imply 
\begin{subequations}\label{6}
\begin{align} \label{6a}
(\alpha_i - \alpha_j) f_{ij}^{(1)} &= f(\alpha_i) - f(\alpha_j),
\\
(\alpha_i - \alpha_j ) (\alpha_j - \alpha_k ) (\alpha_k - \alpha_i) f_{ijk}^{(2)} &=
(\alpha_k - \alpha_j ) f(\alpha_i)  + (\alpha_i - \alpha_k )  f(\alpha_j) 
+(\alpha_j - \alpha_i ) f(\alpha_k), 
\end{align}
\end{subequations}
These may be expressed in the obvious alternate form, 
 \begin{subequations}\label{7}
\begin{align} \label{7a}
 f_{ij}^{(1)} &= \frac{ f(\alpha_i) - f(\alpha_j)} { \alpha_i - \alpha_j }, \quad i\ne j, 
\\
 f_{ijk}^{(2)} &=
 \frac{ f(\alpha_i)}{(\alpha_i - \alpha_j )(\alpha_i - \alpha_k )}
+\frac{ f(\alpha_j)}{(\alpha_j - \alpha_k )(\alpha_j - \alpha_i )}
+\frac{ f(\alpha_k)}{(\alpha_k - \alpha_i )(\alpha_k - \alpha_j )},
\quad i\ne j\ne k \ne i  . 
\end{align}
\end{subequations}
These relations do not, however, provide all the coefficients required to define the derivatives according to \rf{4a} and \rf{4b}.   Specifically, the  coefficients $f_{ii}^{(1)}$, $ f_{ikk}^{(2)}$ and $ f_{kkk}^{(2)}$ are not defined, but  may  be obtained from \rf{6} by taking the limits $\alpha_j \rightarrow \alpha_i$ etc. and using l'Hopital's rule.  The result is $f_{ii}^{(1)} = f'(\alpha_i)$ in agreement with \rf{45}, and 
\begin{equation} \label{050}
 f_{iik}^{(2)} = 
 \begin{cases} 
 \frac{ f'(\alpha_i)}{\alpha_i - \alpha_k } -  
\frac{ f(\alpha_i) - f(\alpha_k)}{(\alpha_i - \alpha_k )^2},  & 
 i\ne k ,  \\
\frac12 f''(\alpha_i) , &   i=k,  
\end{cases}
\end{equation}
in agreement with Theorem \ref{thm1} and with \rf{47} from  \cite{Ogden84}.

How does this generalize to the $n^{th}$ derivative, and is it correct to take the limits as described to find the ``missing" coefficients?  In partial answer to the first question, 
 note that Lemma \ref{lem2} applied  to the function $f$ represented as a  power series, \rf{810}, suggests the equivalence
 \begin{equation} \label{093}
\sum\limits_{i=1}^{n+1} 
 \frac{f(x_i)}{\prod\limits_{\stackrel{j=1}{j\ne i}}^{n+1}  (x_i-x_j)}   
 \quad \rightarrow \quad f_{i_1,i_2,\ldots ,  i_{n+1}}^{(n)} ,\quad n\ge 3, 
 \end{equation} 
where the $x_i$ and $x_j$ are replaced by $\alpha_i$ and $\alpha_j$.  However, no matter what combination of indices are taken for the coefficients on the right hand side, there will be pairs of indices with the same values, e.g. $f_{ijkk}^{(3)}$.  In this sense  the cases 
$n=1$ and $n=2$ are special because it is possible to have indices that are distinct, which enables us to write the identities \rf{7}.  This does not apply for $n\ge 3$.   The indicial overlap must be taken into account in order to arrive at the correct form of \rf{093}.  Not surprisingly, the correct result is equivalent to taking the appropriate limit of the left member of \rf{093}.  Consider the third order derivative,  for instance.  The coefficient 
$f_{ijkk}^{(3)}$ corresponds to $x_1= \alpha_i$, $x_2= \alpha_j$, 
$x_3= \alpha_k$  in the limit as  $x_4\rightarrow \alpha_k$.  That is, 
\begin{align} \label{09}
f_{ijkk}^{(3)} = & \lim_{x\rightarrow \alpha_k} \bigg[
 \frac{ f(\alpha_i)}{(\alpha_i - \alpha_j )(\alpha_i - \alpha_k )(\alpha_i - x )}
+\frac{ f(\alpha_j)}{(\alpha_j - \alpha_k )(\alpha_j - x )(\alpha_j - \alpha_i )}
\nonumber \\
& \qquad 
+\frac{ f(\alpha_k)}{(\alpha_k - \alpha_i )(\alpha_k - \alpha_j )(\alpha_k - x )}
+\frac{ f(x)}{(x - \alpha_i )(x - \alpha_j )(x - \alpha_k )}
\bigg] 
\nonumber \\
= & \frac{1}{(\alpha_i - \alpha_j) }
\bigg[
\frac{ f(\alpha_i)- f(\alpha_k)}{ (\alpha_i - \alpha_k)^{2} } -
\frac{ f(\alpha_j)- f(\alpha_k)}{ (\alpha_j - \alpha_k)^{2} }\bigg]
+    
\frac{ f'(\alpha_k) }{ (\alpha_i - \alpha_k)(\alpha_j - \alpha_k) },
\end{align}
in agreement  with \rf{1d} for $n=3$.  The remaining coefficients $f_{ikkk}^{(3)}$, $f_{iikk}^{(3)}$ and $f_{kkkk}^{(3)}$ can be obtained by taking further limits. 
This method of evaluating the coefficients, while correct, is  particularly tedious as $n$ becomes larger and  more limits must be taken.  

Let us therefore  consider the form of $ f_{i_1,i_2,\ldots ,  i_{n+1}}^{(n)}$.  First, it is clear that 
the coefficient is unchanged under any permutation of the indices  $i_1,i_2,\ldots ,  i_{n+1} $.  In other words, it is totally symmetric in the $n+1$ indices.  Secondly, the indices assume only three different values: $1$, $2$ and $3$. Since  there are $n+1$ indices then apart from the first and second derivatives $f^{(1)}_{ij}$ $(n=1)$  and $f^{(1)}_{ijk}$ $(n=2)$ there will always be coincident indices present in the coefficient $ f_{i_1,i_2,\ldots ,  i_{n+1}}^{(n)}$.  These are indices with the same value whether that value is $1$, $2$ or $3$.  For $n\ge 3$ there are at least $n-1$ coincident indices.   
Suppose the index $k$, distinct from $i$ and $j$, occurs $\nu_k>1$ times in 
$\{i_1,i_2,\ldots ,  i_{n+1}\}$.  Similarly, the indices $i$ and $j$ each occur $\nu_i>0$ and $\nu_j>0$ times.  This suggests the change of notation  
\begin{equation} \label{595}
f_{i_1,i_2,\ldots ,  i_{n+1}}^{(n)} \quad \rightarrow \quad
f_{ijk }^{\nu_i ,\nu_j ,\nu_k }, \qquad \text{with  }\nu_i +\nu_j +\nu_k = n+1. 
\end{equation}
We are now ready to describe  two  approaches for calculating these coefficients.

\subsection{A closed form expression}

 We first note that the coefficient $ f_{i_1,i_2,\ldots ,  i_{n+1}}^{(n)}$ of \rf{4c} can be evaluated as a contour integral.  This follows from \rf{093} and the integral identity  
 \begin{equation} \label{0900}
 \sum\limits_{i=1}^{n+1} 
 \frac{f(x_i)}{\prod\limits_{\stackrel{j=1}{j\ne i}}^{n+1}  (x_i-x_j)}   
 = \frac1{2\pi i} \int\limits_C 
 \frac{\dd z\, f(z) }{ 
 \prod\limits_{l=1}^{n+1}  (z-x_l)},  
 \end{equation} 
 where  $C$ is any contour in the complex plane that encloses the set of points $\{ x_1,x_2, \ldots , x_{n+1}\}$.    The integrand  in \rf{0900} must be modified   to  account for multiple occurrences of the indices, in   the notation of \rf{595}.  We  obtain the integral identity
 \begin{equation} \label{0901}
 f_{ijk }^{\nu_i ,\nu_j ,\nu_k } = 
 \frac1{2\pi i} \int\limits_C 
 \frac{\dd z\, f(z) }{ 
   (z-\alpha_i)^{\nu_i}
   (z-\alpha_j)^{\nu_j}
   (z-\alpha_k)^{\nu_k} },
 \end{equation} 
which can be evaluated by residues, with the final result \rf{305}.

\subsection{Interpolation interpretation}

The following result  is key: 
\begin{lemma}\label{lem3}
If     $\{ x_1,x_2,  \ldots , x_{n+1}\}$ are $n+1$ distinct points then
\begin{equation} \label{443}
\sum\limits_{i=1}^{n+1} \frac{ f(x_i)} {\prod\limits_{\stackrel{k=1}{k\ne i}}^{n+1}  (x_i-x_k)}
= \frac{1}{n!} \frac{\dd^n ~}{\dd x^n} P(x) , 
\end{equation}
where $P(x )$ is the Lagrange interpolating polynomial of the $n+1$ function values $f (x_{i})$, $i=1,2,\ldots, n+1$.
\end{lemma}
The proof follows from the explicit form of the interpolating polynomial of degree $n$, 
\begin{equation} \label{442}
P(x) = 
\sum\limits_{i=1}^{n+1} f (x_i) L_i (x) ,  
\end{equation}
where 
$L_i (x)$, $i=1,2,\ldots, n+1$, are the Lagrange   polynomials  satisfying  $L_i (x_j)= \delta_{ij}$, 
\begin{equation} \label{441}
L_i (x) =  {\prod\limits_{\stackrel{j=1}{j\ne i}}^{n+1}  (x-x_j)} \bigg\slash \,  {\prod\limits_{\stackrel{k=1}{k\ne i}}^{n+1}  (x_i-x_k)} .
\end{equation}
Differentiating $P(x)$ of \rf{442}  yields  the identity \rf{443}. 

This connection between the coefficients and polynomial  interpolation provides an alternative  means to calculate the coefficient $ f_{ijk }^{\nu_i ,\nu_j ,\nu_k } $.  
The interpolation problem that must be solved is not the straightforward one of Lemma \ref{lem3} with a single datum at each of $n+1$ distinct points.
Some  of the points can  have more than one piece of data.  Specifically, 
associated with $x=\alpha_l$ where $l=i$, $j$ and  $k$, will be the $\nu_l$ data comprising $f(\alpha_l)$, $f'(\alpha_l)$, \ldots, $f^{(\nu_l -1)}(\alpha_l)$.  In this way  there is a well defined and unique interpolating polynomial of degree $x^n$ for each $f_{ijk }^{\nu_i ,\nu_j ,\nu_k }$ with 
\begin{equation} \label{59}
f_{ijk }^{\nu_i ,\nu_j ,\nu_k } = \frac{1}{n!} \frac{\dd^n ~}{\dd x^n} P(x ) . 
\end{equation}

  The interpolation problem that must be solved does not fit into the standard category of Lagrange interpolating polynomials, since some or many of the points are coincident.  However, the procedure for generating the  polynomial is not difficult. Appendix \ref{appa} illustrates the method  by example, for the case of $f$ known at three points along with all derivatives  at one of these points up to  $f^{(n-2)}$. 

\subsection{The number of independent coefficients}

The final detail is a count of the number  $N(n)$ of independent expressions  of the form $f_{ijk }^{\nu_i ,\nu_j ,\nu_k } $ required
for the  $n^{th}$ derivative.   This   is the number of elements in the set of  integer triples $\{ I,J,K \in \tens{Z}: \, 0\le I\le J\le K, \, I+J+K=n+1\}$. That is, 
\begin{equation} \label{0383}
N(n) = \dim \, \{ I,J \in \tens{Z}: \, 0\le I\le J\le n+1-I-J\}.
\end{equation}
This defines a region $A(n)$ on  the integer grid of $\{I,J\}$ bounded by the three straight lines $L_1: I=0$,  $L_2: I-J=0$ and $L_3: I+2J=n+1$.  $N(n)$ is the number of points in $A(n)$, including on its boundary.  Enumeration yields 
\begin{equation} \label{433}
N(n) = \lfloor \frac{ (n+4)^2 + 4}{12} \rfloor  
\end{equation}
where the floor function $\lfloor x \rfloor $ is the integral part of $x$. 

This completes the proof of Theorem \ref{thm1}.   Note that the two methods discussed for calculating the coefficients, based on the contour integral and the interpolation problem, are both related to divided differences \cite{Donoghue}.

\section{The inverse gradient}\label{sec4}

\subsection{Definition}
 
The inverse tensor function 
  $ \nabla^{-1} f ({\bd A}) \in \tens{L}$in  is defined by 
  $(\nabla f ) \nabla^{-1} f = (\nabla^{-1} f ) \nabla f = \tens{I}$. 
As an example application  let $f^\text{inv}$ be the inverse of $f$, such that 
$f\circ f^\text{inv} ({\bd A}) = f^\text{inv}\circ f  ({\bd A})= {\bd A}$.  It follows from \rf{557} that 
\begin{equation} \label{556}
\nabla  f^\text{inv} ({\bd A})  = 
\nabla^{-1}  f\big( f^\text{inv}({\bd A})\big)  . 
\end{equation}
It is clear from \rf{45} that the definition of $ \nabla^{-1} f ({\bd A})$ is problematic if $f'(\alpha_i)$ vanishes. This possibility is precluded  by restricting consideration to strictly monotonic \emph{strain measure functions} acting on positive definite tensors.

\subsection{Application to strain measure functions}
 
 The function $f$ is a  \emph{strain measure} \cite{Hill78,Scheidler91} if it is a smooth function 
$ f : \tens{R}^+ \rightarrow \tens{R}$ which satisfies 
\begin{equation} \label{004}
f(1)=0,\qquad f'(1)=1, \qquad f'>0. 
\end{equation}
For the remainder of the paper we restrict attention to strain measure functions. 
Examples will be presented for the  Seth-Hill strain measure functions, 
 \begin{equation} \label{777}
 f^{[m]}(x) =  m^{-1} (x^m-1).  
 \end{equation}
The inverse $f^\text{inv}$ is well defined since strain measure functions are one-to-one.  Hence, the inverse gradient is relevant to calculating the gradient of the inverse function, via the identity \rf{557}. 

Strain measure functions act on positive definite strain tensors ${\bd S}\in P$sym, with 
\begin{equation} \label{660}
{\bd S} = \sum\limits_{i=1}^d\lambda_i{\bd S}_i,
\qquad
{\bd I} = \sum\limits_{i=1}^d{\bd S}_i,
\qquad  \lambda_i>0, \qquad {\bd S}_i{\bd S}_j= {\bd S}_j{\bd S}_i = \delta_{ij}{\bd S}_i ,
\end{equation}
where $d\le 3$ is now the eigen-index of ${\bd S}$. 

Accordingly, we extend the concept of positive definiteness to fourth order tensors.  It is necessary to  define an  inner product on Lin, which is done in the usual manner as 
${\bd A}\cdot {\bd B} = \tr({\bd A} {\bd B}^t)$.
Then ${\tens L}$in is the space of fourth order tensors acting on Lin,  
 ${\bd X} \rightarrow {\tens A}{\bd X}$ with  transpose ${\tens A}^t$ such that 
 ${\bd Y}\cdot {\tens A}{\bd X} = {\bd X}\cdot{\tens A}^t{\bd Y}$ for all
 ${\bd X}$, ${\bd Y}\in$Lin.   
 The vector space may be decomposed ${\tens L}$in $= {\tens S}$ym$ \oplus  {\tens S}$kw 
where  ${\tens S}$ym and ${\tens S}$kw denote the spaces of symmetric (${\tens A}^t={\tens A}$) and skew-symmetric (${\tens A}^t=-{\tens A}$) tensors, respectively.   
Any ${\tens A}\in  {\tens L}$in can be uniquely partitioned into symmetric and skew parts: 
${\tens A}= {\tens A}^{(+)}+{\tens A}^{(-)}$, where 
${\tens A}^{(\pm)} = ({\tens A}\pm {\tens A}^t)/2$. 
 The identity ${\tens I}$ satisfies ${\tens I}{\bd X} = {\bd X} $ for all ${\bd X}\in$ Lin.   The product ${\tens A} {\tens B}\in \tens{L}$in is defined by 
${\bd Y}\cdot{\tens A} {\tens B}{\bd X} =   ({\tens A}^t{\bd Y}) \cdot {\tens B}{\bd X}$. 
${\tens P}$sym  is the space of positive definite fourth order tensors:  ${\tens A}\in {\tens P}$sym  iff ${\bd X}\cdot {\tens A}{\bd X} >0$, for all nonzero $
{\bd X} \in$ Sym.  

The spectral form of ${\tens A}\in {\tens P}$sym is 
\begin{equation} \label{-99}
{\tens A} = \sum\limits_{I=1}^D a_I {\tens A}_I, 
\qquad  
{\tens I} = \sum\limits_{I=1}^D {\tens A}_I, 
\qquad  
{\tens A}_I{\tens A}_J = \begin{cases}
{\tens A}_I & I=J,\\
0 , & I\ne J,  
\end{cases}
\end{equation}
where $D$ is the eigen-index, $a_I>0$ and ${\tens A}_I\in \tens{P}$sym. 
Functions of fourth order tensors  can therefore be defined according to 
\begin{equation} \label{845}
f( {\tens A}) =   \sum\limits_{I=1}^D f(a_I) {\tens A}_I . 
\end{equation}

For any ${\bd A}\in$Psym and its eigen-tensors ${\bd A}_i$, $i=1,\ldots , d$ of \rf{-71},  define the associated set of $D=\frac12 d(d+1)$ basis tensors for $\tens{S}$ym by 
\begin{equation} \label{038}
 \tens{A}_I =  \begin{cases} 
    \bd{A}_I\boxtimes \bd{A}_I, &  I=1,\ldots , d, 
  \\
   \bd{A}_i\boxtimes \bd{A}_j+\bd{A}_j\boxtimes \bd{A}_i ,
 & I=d+1 , \ldots , D =\frac12 d(d+1). 
  \end{cases} 
  \end{equation}
Here $I=d+1 , \ldots , D$ corresponds to  distinct pairs $(i,j)$ with $i< j$. 
Thus, $D=6$ for $d=3$,  while $D=3$ for $d=2$ and  $D=1$ in the trivial case of  $d=1$. 
It may be readily checked that $\tens{A}_I \tens{A}_J = \tens{A}_I \delta_{IJ}$.  
The identity ${\tens I}   =  {\bd I} \boxtimes {\bd I}$  implies 
the partition of unity 
${\tens I} 
 =
 \sum\limits_{I=1}^D  
  \tens{A}_I$.

In particular, 
\begin{lemma}\label{lem3a}
The gradient of a strain measure function and its inverse are positive definite fourth order tensors, i.e. 
$\nabla f ({\bd A}),\, \nabla^{-1} f ({\bd A}) \in \tens{P}$sym, and 
\begin{equation} \label{48}
(\nabla f ) \nabla^{-1} f = (\nabla^{-1} f ) \nabla f = \tens{I}. 
 \end{equation} 
 The spectral forms are 
 \begin{subequations}\label{362}
 \begin{equation} \label{362a}
 \nabla f ({\bd A}) = \sum\limits_{I=1}^D f_I \tens{A}_I, 
 \qquad 
 \nabla^{-1} f ({\bd A}) = \sum\limits_{I=1}^D f_I^{-1} \tens{A}_I, 
 \end{equation}
 with 
 \begin{equation} \label{372}
 f_I ({\bd A}) = \begin{cases} 
 f'(\alpha_I),  &  I=1,\ldots , d,  \\ 
 \frac{ f(\alpha_i) - f(\alpha_j)}{ \alpha_i - \alpha_j}, & I=d+1 , \ldots , D. 
 \end{cases} 
 \end{equation}
 \end{subequations}
 \end{lemma}
The proof is evident  from   the spectral decompositions.
 The positivity of  
 $ f_I$, $I=1, \ldots , D$ 
is  a  consequence of $\alpha_i >0$, $i=1,\ldots , d$ combined with the monotonicity of $f$ and the mean value theorem.   
 The positive definite nature of $ \nabla f ({\bd A})$ implies that it has a unique positive definite inverse.

 The  gradient function and its inverse has an  interesting alternative representation for integer and fractional values of $m$, respectively: 
 \begin{lemma}\label{lem4}
 For integer values of $m\ne 0$, 
 \begin{subequations}
 \begin{align} \label{368}
 \nabla f^{[m]} ({\bd A}) &= \frac{1}{|m|} \begin{cases}
 \sum\limits_{k=1}^{m}
 \bd{A}^{m-k}\bxt \bd{A}^{k-1}, & m >0,
 \\
 \sum\limits_{k=m+1}^{0}
 \bd{A}^{m-k}\bxt \bd{A}^{k-1}, & m <0.
 \end{cases}
 \\
 \nabla^{-1} f^{(\frac1{m})} ({\bd A}) &= \frac{1}{|m|} \begin{cases}
 \sum\limits_{k=1}^{m}
 \bd{A}^{1-\frac{k}{m}}\bxt \bd{A}^{\frac{k}{m}- \frac{1}{m}}, & m >0,
 \\
 \sum\limits_{k=m+1}^{0}
 \bd{A}^{1-\frac{k}{m}}\bxt \bd{A}^{\frac{k}{m}- \frac{1}{m}}, & m <0.
 \end{cases}
 \label{368b}
 \end{align}
 \end{subequations}
 \end{lemma}
The proof uses the expansion of eq. \rf{362} with $f_I \rightarrow f^{[m]}_I$
 where
 \begin{equation} \label{366}
 f^{[m]}_I ({\bd A}) = \begin{cases} 
 \alpha_I^{m-1}, &  I=1,\ldots , d   \\ 
 \frac{1}{m}\big( \frac{ \alpha_i^m - \alpha_j^m}{ \alpha_i - \alpha_j}\big),
 & I=d+1, \ldots , D.
 \end{cases} 
 \end{equation}
  Consider the case of positive $m$.  Then using the spectral form for  $\bd A$  
 \begin{align} \label{369}
\frac{1}{m}
 \sum\limits_{k=1}^{m}
 \bd{A}^{m-k}\bxt \bd{A}^{k-1}
 &= \frac{1}{m} \sum\limits_{i,j=1}^d  
 {\bd A}_i \bxt {\bd A}_j  \sum\limits_{k=1}^{m}
 \alpha_i^{m-k} \alpha_j^{k-1}
 \nonumber \\
  &=  \sum\limits_{i=1}^d  \alpha_i^{m-1}{\bd A}_i \bxt {\bd A}_i
  + \frac{1}{m}\sum\limits_{\stackrel{ i,j=1}{i\ne j}}^d  
 {\bd A}_i \bxt {\bd A}_j \alpha_i^{m}\alpha_j^{-1}
 \sum\limits_{k=1}^{m} \big( \frac{\alpha_j}{\alpha_i}\big)^k 
 \nonumber \\
  &=  \sum\limits_{I=1}^d  f^{[m]}_I \tens{A}_I
  + \frac{1}{m}\sum\limits_{\stackrel{ i,j=1}{i\ne j}}^d  
 {\bd A}_i \bxt {\bd A}_j  \big( \frac{\alpha_i^m -\alpha_j^m}{\alpha_i-\alpha_j}\big) 
 \nonumber \\
  &=  \nabla f^{[m]} ({\bd A}), 
 \end{align}
 where the identity $x+ x^2 + \ldots + x^m=x(1-x^m)/(1-x)$ has been used.  The proof for negative $m$ is similar.  The results for $\nabla^{-1} f^{(\frac1{m})} ({\bd A})$ are a consequence of the identity
 \begin{equation} \label{32}
 f^{(\frac1{m})}_I ({\bd A}) = 1/f^{[m]}_I ({\bd A}^\frac1{m}),
 \qquad m\ne 0.
 \end{equation}
 
 Some examples of interest, first for integer values of $m$:
 \begin{subequations}
  \begin{align} \label{379}
   \nabla f^{(1)} ({\bd A}) &= \bd{I} \bxt \bd{I}  =\tens{I},
   \\
    \nabla f^{(-1)} ({\bd A}) &= \bd{A}^{-1}\bxt \bd{A}^{-1},
   \\
    \nabla f^{(2)} ({\bd A}) &= \frac12 ( \bd{I}\bxt \bd{A}  + \bd{A}\bxt \bd{I}), \label{379c}
   \\
   \nabla f^{(-2)} ({\bd A}) &= \frac12 ( \bd{A}^{-2}\bxt \bd{A}^{-1}  + \bd{A}^{-1}\bxt \bd{A}^{-2}).
   \end{align}
 \end{subequations}
 The inverse tensors 
 \begin{subequations}
  \begin{align} \label{3670}
   \nabla^{-1} f^{(1)} ({\bd A}) & =\tens{I},
   \\
    \nabla^{-1} f^{(-1)} ({\bd A}) &= \bd{A} \bxt \bd{A} ,
   \end{align}
 \end{subequations}
 follow by observation, as do the general relations 
 \begin{subequations}
  \begin{align} \label{371}
   \nabla f^{(-m)} ({\bd A}) &= 
    \bd{A}^{-m}\bxt \bd{A}^{-m} \, \nabla f^{[m]} ({\bd A}) ,
   \\
    \nabla^{-1} f^{(-m)} ({\bd A}) &= 
     \bd{A}^{m}\bxt \bd{A}^{m} \, \nabla^{-1} f^{[m]} ({\bd A}) ,
  \end{align}
 \end{subequations}
 which express tensors in terms of their counterparts of opposite sign. 
 
${\bd A}{\bd X}+{\bd X}{\bd A}={\bd C}$  is an important and common  equation in mechanics \cite{Scheidler94,Rosati00}.  The  solution may be written simply and succinctly 
as ${\bd X} = ({\bd A}\bxt{\bd I} + {\bd I} \bxt{\bd A} )^{-1}{\bd C}$
where  \cite{Jog06} (eq. (22))
\begin{equation} \label{68}
({\bd A}\bxt{\bd I} + {\bd I} \bxt{\bd A} )^{-1}
 = \sum\limits_{i,j=1}^d \, 
 (\alpha_i + \alpha_j)^{-1} {\bd A}_i\bxt{\bd A}_j .
\end{equation}
  In the current notation  this may be expressed 
${\bd X} = \frac12 \nabla^{-1} f^{(2)}({\bd A} ){\bd C}$, see eq. \rf{379c}. 
In fact, the identity \rf{368} implies that  the solution to 
\begin{equation} \label{381}
{\bd A}^{m-1}{\bd X} + {\bd A}^{m-2}{\bd X}{\bd A} +{\bd A}^{m-3}{\bd X}{\bd A}^2 +
\ldots + {\bd X}{\bd A}^{m-1} = {\bd C}, 
\end{equation}
is
\begin{equation} \label{-38}
{\bd X} = \frac1{m} \nabla^{-1} f^{[m]}({\bd A} )\, {\bd C} = 
\sum\limits_{ i,j=1}^d  
 \big( \frac{\alpha_i-\alpha_j}{\alpha_i^m -\alpha_j^m} \big) 
 {\bd A}_i \bxt {\bd A}_j \, {\bd C}. 
\end{equation}
The solution to ${\bd A}{\bd X}+{\bd X}{\bd A}={\bd C}$ can be written in a form that does not require the spectral split, e.g.  
${\bd X} = [ 2(I_1 I_2 - I_3)]^{-1}\, 
\big[ 
(I_1^2 - I_2){\bd C}  + {\bd A}{\bd C}{\bd A} - ({\bd A}^2 {\bd C} + {\bd C}{\bd A}^2)  + I_1I_3{\bd A}^{-1}{\bd C}{\bd A}^{-1}
\big],
$ \cite{Rosati00} (eq. (27))
where $I_1, I_2, I_3$ are the invariants of ${\bd A}$.   Similar expressions, although far more complicated,  can be generated for the solution to \rf{381} and related equations, but we leave such matters aside. 

Examples of fractional powers include the well known case of $\frac12$ \cite{Rosati00}, 
for which \rf{368b} gives
\begin{equation} \label{-90}
  \nabla^{-1} f^{(\frac12 )} ({\bd A}) = \frac12( \sqrt{\bd{A}} \bxt \bd{I} +\bd{I}\bxt\sqrt{\bd{A}})
  \quad
  \Rightarrow \quad
  \nabla  \sqrt{\bd{A}} = ( \sqrt{\bd{A}} \bxt \bd{I} +\bd{I}\bxt\sqrt{\bd{A}})^{-1}.
\end{equation}
Similarly, 
\begin{equation} \label{191}
  \nabla^{-1} f^{(-\frac12 )} ({\bd A}) = \frac12( \sqrt{\bd{A}} \bxt \bd{A} +\bd{A}\bxt\sqrt{\bd{A}})
  \quad
  \Rightarrow \quad
  \nabla {\bd{A}}^{-1/2}  = -( \sqrt{\bd{A}} \bxt \bd{A} +\bd{A}\bxt\sqrt{\bd{A}})^{-1},
\end{equation}
and
\begin{subequations}
  \begin{align} \label{-93}
   \nabla^{-1} f^{(\frac13)} ({\bd A}) &= 
  \frac13(  \bd{A}^{2/3} \bxt \bd{I}+\bd{A}^{1/3} \bxt \bd{A}^{1/3}
  +\bd{I} \bxt \bd{A}^{2/3}),
  \\
  \nabla^{-1} f^{(-\frac13)} ({\bd A}) &= 
  \frac13(  \bd{A}^{1/3} \bxt \bd{A}+\bd{A}^{2/3} \bxt \bd{A}^{2/3}
  +\bd{A} \bxt \bd{A}^{1/3}),
  \end{align}
 \end{subequations}
 imply expressions for $\nabla {\bd{A}}^{\pm 1/3}$, etc.

 Note that the case of $m=0$, corresponding to 
 \begin{equation} \label{460}
\nabla f^{(0)} ({\bd A}) = \nabla \ln ({\bd A}) , 
 \end{equation}
 is specifically excluded in Lemma \ref{lem3}.  Although there is no expression for 
 $\nabla \ln ({\bd A})$ analogous to \rf{368},  its inverse can be  represented as follows:
 \begin{lemma}\label{lem5}
 The tensor  $ \nabla^{-1} \ln ({\bd A}) = \nabla^{-1} f^{(0)} ({\bd A}) $, which is the inverse of  $\nabla \ln ({\bd A})$, is
 \begin{equation} \label{33}
\nabla^{-1} \ln ({\bd A}) = 
\int\limits_0^1 \dd x\, {\bd A}^x \bxt {\bd A}^{1 - x}
, \qquad \text{for any }{\bd A}\in \text{Psym} . 
 \end{equation} 
 \end{lemma}
 The proof follows from  \rf{368b} by taking the limit of $m\rightarrow \infty$.  An alternative and direct proof can be obtained using the spectral form of ${\bd A}$, 
\begin{align} \label{34}
\int\limits_0^1 \dd x\, {\bd A}^x \bxt {\bd A}^{1 - x} 
&= \sum\limits_{i,j=1}^d 
 {\bd A}_i \bxt {\bd A}_j\, \int\limits_0^1 \dd x\, \alpha_i^x\alpha_j^{1-x}
\nonumber \\
&= 
 \sum\limits_{i,j=1}^d \, \frac{ \alpha_i - \alpha_j}{ \ln\alpha_i - \ln\alpha_j}
 {\bd A}_i \bxt {\bd A}_j, 
\end{align}
where the ratio becomes $\alpha_i$ if $i=j$.    The right member of \rf{34} is obviously the inverse of  $\nabla \ln ({\bd A})$ from eq. \rf{45}.

 \subsection{The tensors $\tens {J}$,  $\tens {J}^*$,  $\tens {K}$, and the solution to $AX-XA=Y$}
  
Consider the equation 
\begin{equation} \label{60}
{\bd A}{\bd X} - {\bd X} {\bd A} = {\bd Y},\qquad {\bd A}\in \text{Psym} , 
\end{equation}
for the unknown ${\bd X}$ in terms of  ${\bd Y}$ which  is either symmetric or skew, and ${\bd X}$ is of the opposite parity      
\cite{Dui06}.
The equation can be written 
\begin{equation} \label{61}
\tens{J}({\bd A}){\bd X}  = {\bd Y},
\end{equation}
where  
\begin{equation} \label{62}
\tens{J}({\bd A})  \equiv {\bd A}\bxt{\bd I} - {\bd I} \bxt{\bd A}, \qquad {\bd A}\in \text{Psym} .
\end{equation}

We will only consider $\tens{J}({\bd A})$ for symmetric $\bd A$, implying 
$\tens{J}\in \tens{S}$ym and $\tens{J}$ maps Sym$\rightarrow$Skw and Skw$\rightarrow$Sym. Therefore, 
$\tens{J}$ does not possess eigenvalues,  eigenvectors or an inverse in the usual sense.  
 However, it is possible to define the pseudo-inverse, or equivalently the Moore-Penrose inverse,  $\tens{J}^*$, which satisfies
\begin{align} \label{631}
\tens{J}\tens{J}^*\tens{J} & = \tens{J},
\nonumber \\
\tens{J}^*\tens{J} \tens{J}^* & = \tens{J}^* . 
\end{align}
Further understanding comes from consideration of the spectral form 
\begin{equation} \label{63}
\tens{J}({\bd A})  = \sum\limits_{i,j=1}^d  (\alpha_i - \alpha_j) {\bd A}_i\bxt{\bd A}_j .
\end{equation}
This representation implies that ${\bd A}_i$ are null vectors, i.e. 
$\tens{J}({\bd A}){\bd A}_i = 0$, $i=1,\ldots , d$. 
Noting that 
\begin{equation} \label{64}
\tens{J}({\bd A})  [{\bd A}_i\bxt{\bd A}_j \pm {\bd A}_j\bxt{\bd A}_i]
=  (\alpha_i - \alpha_j) [{\bd A}_i\bxt{\bd A}_j \mp {\bd A}_j\bxt{\bd A}_i] , 
\end{equation}
gives
\begin{equation} \label{65}
\tens{J}^*({\bd A})  [{\bd A}_i\bxt{\bd A}_j \pm {\bd A}_j\bxt{\bd A}_i]
=  (\alpha_i - \alpha_j)^{-1} [{\bd A}_i\bxt{\bd A}_j \mp {\bd A}_j\bxt{\bd A}_i] , \quad \text{for } i\ne j \text{ only}.
\end{equation}
The caveat $i\ne j$ is crucial, and indicates that the non-null parts of $\tens{J}$ and its pseudo-inverse define  maps between three dimensional subspaces of Sym and Skw.  In particular, 
\begin{equation} \label{66}
\tens{J}^*({\bd A})  = 
({\bd A}\bxt{\bd I} - {\bd I} \bxt{\bd A} )^* = \sum\limits_{\substack{i,j=1 \\ i\ne j}}^d \, 
 (\alpha_i - \alpha_j)^{-1} {\bd A}_i\bxt{\bd A}_j .
\end{equation}
This formula for the pseudo-inverse clearly satisfies \rf{631}.  
 It is interesting to compare the form of \rf{66} with \rf{68}. 

Note that 
\begin{equation} \label{67}
\tens{J}^*\tens{J} = \tens{J}\tens{J}^* = 
\tens{I} - \sum\limits_{I =1}^d   {\tens A}_I  = \sum\limits_{I =d+1}^D   {\tens A}_I. 
\end{equation}
Define $\tens{K} ({\bd A})$, ${\bd A}\in$Psym,  and its pseudoinverse,
\begin{equation} \label{674}
\tens{K} ({\bd A})= \sum\limits_{i =1}^d   \alpha_i {\tens A}_i,
\qquad 
\tens{K}^* ({\bd A})= \sum\limits_{i =1}^d   \alpha_i^{-1} {\tens A}_i = \tens{K} ({\bd A}^{-1}),
\end{equation}
then 
\begin{equation} \label{6741}
\tens{J}({\bd A})\tens{J}^*({\bd A}) +  \tens{K}({\bd A})\tens{K}^*({\bd A}) = 
\tens{I} \qquad \text{for any } {\bd A} \in \text{Psym}. 
\end{equation}
The tensors $\tens{J}$ and $\tens{K}$ allow us to  write  $\nabla f({\bd A})$ and its inverse in a succinct manner as, respectively,  
\begin{subequations}
\begin{align} \label{673}
\nabla f({\bd A}) &=  \tens{K} \big(f'({\bd A}) \big)
+ \tens{J}^*({\bd A})  \tens{J}\big( f({\bd A})\big) ,
\\
 \nabla^{-1} f({\bd A}) &= \tens{K}^* \big(f'({\bd A}) \big)
+ \tens{J}({\bd A})  \tens{J}^*\big( f({\bd A})\big). 
\end{align}
\end{subequations}

\appendix
\Appendix{Solving the  interpolation problems}\label{appa}

Consider three points, say $x_1$,  $x_2$ and  $x_3$. 
Given $f(x_1)$, $f(x_2)$ and $f^{(l)}(x_3)$, $l=0,1,2,\ldots , n-2$ find the interpolating polynomial $P(x)$ of degree $n$.   For simplicity, but no lack of generality, take $x_3 = 0$, and consider the ansatz
\begin{equation} \label{033}
P(x) = \sum\limits_{l=0}^{n} p_l x^l, \qquad p_l= \frac{1}{l!} f^{(l)}(0), \quad l=1,2,\ldots , n-2.
\end{equation}
This satisfies all $n-1$  conditions at $x=0$, and the conditions at  $x_1$ and   $x_2$ are met if $P(x_1) = f(x_1)$ and $P(x_2) = f(x_2)$. These imply a pair of    linear equations in 
$p_{n-1}$ and $p_n$ which are easily solved.  Only the coefficient $p_n$ is required as this  determines the $n^{th}$ derivative of $P(x)$, 
\begin{equation} \label{-75}
p_n = \frac{1}{x_1-x_2}
\bigg[ \frac{f(x_1)}{x_1^{n-1}}  - \frac{f(x_2)}{x_2^{n-1}} 
- \sum\limits_{l=0}^{n-2} \frac{f^{(l)}(0)}{l!}  \bigg( 
\frac{x_1^l}{x_1^{n-1}} - \frac{x_2^l}{x_2^{n-1}} \bigg)
\bigg].
\end{equation}
This implies eq. \rf{1d}. 

Consider the case where  the values of the function and its first $n-3$ derivatives  are given at one point, say $x=0$, along with $f(x_1)$, $f(x_2)$ and $f'(x_2)$.  Using the same ansatz as \rf{033}  with the three unknowns $p_{n-2}$, $p_{n-1}$ and $p_n$, implies three simultaneous equations.   The solution for $p_n$ is 
\begin{align} \label{-76}
p_n =&  
\frac{1}{(x_1-x_2)^2}
\bigg[ \frac{f(x_1)}{x_1^{n-2}}  - \frac{f(x_2)}{x_2^{n-2}} [ 1- (n-2)(\frac{x_1}{x_2}-1)]
+  \frac{f'(x_2)}{x_2^{n-2}} (x_2-x_1)
\nonumber \\
& \qquad 
- \sum\limits_{l=0}^{n-3} \frac{f^{(l)}(0)}{l!}  \bigg( 
\frac{x_1^l}{x_1^{n-2}} - \frac{x_2^l}{x_2^{n-2}} [ 1- (n-2-l)(\frac{x_1}{x_2}-1)]
\bigg)
\bigg].
\end{align}
This yields the expression  \rf{1e}.

\section*{Acknowledgment}
Thanks to Doron Zeilberger for being true to his motto: ``Who you gonna call?".


\begin{thebibliography}{10}

\bibitem{Carlson86}
D.~E. Carlson and A.~Hoger, \emph{{The derivative of a tensor-valued function
  of a tensor}}, Q. Appl. Math. \textbf{44} (1986), 409--423.

\bibitem{Chadwick71}
P.~Chadwick and R.~W. Ogden, \emph{A theorem of tensor calculus and its
  application to isotropic elasticity}, Arch. Rat. Mech. Anal. \textbf{44}
  (1971), no.~1, 54--68.

\bibitem{Chen04}
Y-C. Chen and G-S. Dui, \emph{The derivative of isotropic tensor functions,
  elastic moduli and stress rate: I. {E}igenvalue formulation}, Math. Mech.
  Solids \textbf{9} (2004), no.~5, 493--511.

\bibitem{Donoghue}
W.~F. Donoghue, Jr., \emph{Monotone matrix functions and analytic
  continuation}, Springer-Verlag, New York, 1974.

\bibitem{Dui06}
Guan-Suo Dui, \emph{Some basis-free formulae for the time rate and conjugate
  stress of logarithmic strain tensor}, J. Elasticity \textbf{83} (2006),
  no.~2, 113--151.

\bibitem{Hill78}
R.~Hill, \emph{{Aspects of invariance in solid mechanics}}, Adv. Appl. Mech.
  \textbf{18} (1978), 1--75.

\bibitem{Itskov02b}
M.~Itskov and N.~Aksel, \emph{A closed-form representation for the derivative
  of non-symmetric tensor power series}, Int. J. Solids Struct. \textbf{39}
  (2002), no.~24, 5963--5978.

\bibitem{Jog06}
C.~S. Jog, \emph{Derivatives of the stretch, rotation and exponential tensors
  in n -dimensional vector spaces}, J. Elasticity \textbf{82} (2006), no.~2,
  175--192.

\bibitem{Kintzel06}
O.~Kintzel and Y.~Bascedilar, \emph{Fourth-order tensors - tensor
  differentiation with applications to continuum mechanics. {I}: Classical
  tensor analysis}, Z. Angew. Math. Mech. \textbf{86} (2006), no.~4, 291--311.

\bibitem{Ogden84}
R.~W. Ogden, \emph{Non-linear elastic deformations}, Ellis Horwood, 1984.

\bibitem{Rosati00}
L.~Rosati, \emph{A novel approach to the solution of the tensor equation
  {AX+XA=H}}, Int. J. Solids Struct. \textbf{37} (2000), no.~25, 3457--3477.

\bibitem{Scheidler91}
M.~Scheidler, \emph{{Time rates of generalized strain tensors. Part I:
  Component formulas}}, Mech. Materials \textbf{11} (1991), 199--210.

\bibitem{Scheidler94}
\bysame, \emph{The tensor equation {AX+XA=F(A,H)}, with applications to
  kinematics of continua}, J. Elasticity \textbf{36} (1994), no.~2, 117--153.

\bibitem{Xiao98}
H.~Xiao, O.~T. Bruhns, and A.~Meyers, \emph{Strain rates and material spins},
  J. Elasticity \textbf{52} (1998), no.~1, 1--41.

\end{thebibliography}

\providecommand{\bysame}{\leavevmode\hbox to3em{\hrulefill}\thinspace}
\providecommand{\MR}{\relax\ifhmode\unskip\space\fi MR }
\providecommand{\MRhref}[2]{%
  \href{http://www.ams.org/mathscinet-getitem?mr=#1}{#2}
}
\providecommand{\href}[2]{#2}

\end{document}